\documentclass[a4paper,12pt,twoside]{article}

\usepackage{microtype}
\usepackage{color}
\usepackage{amssymb,amsmath,amssymb,amsthm}
\usepackage{pdfpages}

\newtheorem{theorem}{Theorem}

\newtheorem{remark}{Remark}
\renewcommand\marginpar[1]{}
\newcommand\comment[1]{}

\definecolor{violet}{rgb}{0.85,0.05,0.85}

\definecolor{dgreen}{rgb}{0,0.6,0.1}

\definecolor{dred}{rgb}{0.8,0,0}

\definecolor{dmagenta}{rgb}{0.6,0,0.6}

\definecolor{dblue}{rgb}{0,0,0.7}

\definecolor{dbrown}{rgb}{0.5,0.25,0}

\definecolor{dorange}{rgb}{1,0.4,0}

\definecolor{dyellow}{rgb}{1,1,0}

\usepackage{cite}

\def\pto{.}
\def\RR{{\mathbb{R}}}
\def\NN{{\mathbb{N}}}

\def\multimathop #1 {\def\arg{#1}%
  \ifx\arg\pto \let\next\relax
  \else
  \def\next{\expandafter
    \def\csname #1\endcsname{\mathop{\rm #1}\nolimits}%
    \multimathop}%
  \fi \next}

\multimathop
dist div dom meas sign supp .

\newenvironment{ele}{}{}
\newcommand{\elena}{\begin{ele}}
\newcommand{\finele}{\end{ele}}
\newcommand{\giu}{\color{blue}}
\newcommand{\finegiu}{\color{black}}
\renewcommand{\giu}{}
\renewcommand{\finegiu}{}
\def\elefin #1{{\color{dgreen}#1}}
\def\pier #1{{\color{dred}#1}}
\def\gius #1{{\color{blue}#1}}
\def\pcol #1{{\color{dred}#1}}

\def\pier #1{#1}
\def\gius #1{#1}
\def\elefin #1{#1}
\def\pcol #1{#1}
\title{A non-smooth regularization\\ of a forward-backward parabolic equation}

\date{}

\begin{document}

\pier{%
\author{}
\thispagestyle{empty}
\begin{center}
\Large{\bf A non-smooth regularization\\ of a forward-backward parabolic equation\footnote{\pier{{\bf Acknowledgment.}\quad\rm
\giu EB and PC acknowledge partial
\finegiu financial support from the MIUR-PRIN Grant 2010A2TFX2 ``Calculus of Variations'', the GNAMPA (Gruppo Nazionale per l'Analisi Matematica, la Probabilit\`a e le loro Applicazioni) of INdAM (Istituto Nazionale di Alta Matematica) and the IMATI -- C.N.R. Pavia.}
\giu GT acknowledges partial financial support from the Italian GNFN (Gruppo Nazionale per la Fisica Matematica) of INdAM through the funding scheme ``Progetto Giovani''.}\finegiu}
\end{center}
\begin{center}
{\large\bf Elena Bonetti$^{(1)}$}\\
{\normalsize e-mail: {\tt elena.bonetti@unipv.it}}\\[.2cm]
{\large\bf Pierluigi Colli$^{(1)}$}\\
{\normalsize e-mail: {\tt pierluigi.colli@unipv.it}}\\[.2cm]
{\large\bf Giuseppe Tomassetti$^{(2)}$}\\
{\normalsize e-mail: {\tt tomassetti@ing.uniroma2.it}}\\[.4cm]
$^{(1)}$
{\small Dipartimento di Matematica ``F. Casorati'', Universit\`a di Pavia}\\
{\small Via Ferrata 1, 27100 Pavia, Italy}\\[.2cm]
$^{(2)}$
{\small Dipartimento di Ingegneria Civile, Universit\`a di Roma ``Tor Vergata''}\\
{\small Via di Tor Vergata 110, 00133 Roma, Italy}\\[.2cm]
\end{center}
}%

\begin{abstract}


\pier{In this paper we introduce a model describing diffusion of species by a suitable 
regularization of a ``forward-backward'' parabolic equation. In particular, we prove 
existence and uniqueness of solutions, as well as continuous dependence on data, for a 
system of partial differential \elefin{\pcol{equations and inclusion}, which may be 
interpreted, e.g., as evolving equation for physical quantities \pcol{such} as}
concentration and chemical potential.
The model deals with a constant mobility
and \elefin{it is recovered from} a possibly non-convex free\pier{-}energy density. \elefin{In particular, we render  a general viscous regularization} via a maximal monotone graph acting on the time derivative of the concentration and presenting a strong coerciveness property.
\\[2mm]
\\[2mm]
}
\medskip
\end{abstract}
\pagestyle{myheadings}
\newcommand\testopari{\sc \pier{Bonetti \ --- \ Colli \ --- \ Tomassetti}}
\newcommand\testodispari{\sc \pier{Regularization of a forward-backward parabolic equation}}
\markboth{\testodispari}{\testopari}
\section{Introduction}
\setcounter{equation}{0}
\elena 
The model we are introducing may be applied to different situations, dealing with diffusion of 
different species (located in some domain $\Omega\subseteq\mathbb{R}^n$ \pier{and}) described in 
terms of concentration. In the following, moving from the classical approach in thermodynamics 
which leads to the well-known Cahn\pier{--}Hilliard equation, we introduce our point of view and make some 
comments on its thermodynamical consistency. In particular, we are focusing on diffusion in solids  
\pier{and we have in mind,} as a possible final application, hydrogen storage in metals.\finele 
\medskip

\noindent\textbf{\elefin{Models for species diffusion}}.  In the simplest setting, a mathematical model describing species diffusion in a solid can be obtained by combining the following three ingredients:

\noindent 1) the \emph{mass-balance law} \pier{for the \emph{concentration} $u$}
\begin{equation}\label{eq:11}
  \begin{aligned}
  &\dot u + \div{\mathbf h}=0;
\end{aligned}
\end{equation}
2) a linear \emph{constitutive prescription} relating  the \emph{flux of diffusant} $\mathbf h$ to the gradient of \emph{chemical-potential} $\mu$ through a \emph{mobility constant} $m$:
\begin{equation}\label{eq:5}
    \mathbf h=-m\nabla \mu,\qquad m>0;
\end{equation}
3) a possibly nonlinear relation between \pier{$\mu$ and $u$}:
\begin{equation}\label{eq:56}
  \mu=\psi'(u),
\end{equation}
dictated by the derivative $\psi'$ of \pier{a} \emph{coarse-grain free energy} $\psi$. 

When diffusion is accompanied by \emph{phase separation}, \elena in general \finele  $\psi$ is \elena assumed \finele to be a non-convex function, a typical choice being the \emph{double-well polynomial potential}:
\begin{equation}\label{eq:50}
  \psi(u)=k(u-1)^2u^2,\qquad k>0.
\end{equation}
\elefin{This kind of specifications} renders
the system \eqref{eq:11}--\eqref{eq:56} \emph{backward parabolic}
\elefin{(e.g., for \eqref{eq:50}, in the region where $0<u<1$)}, an undesirable feature which calls for a suitable regularization.

In order to make the system \elefin{well-posed}, the most popular approach is the so\pier{-}called \emph{elliptic regularization} of \eqref{eq:56}:
\begin{equation}\label{eq:58}
\mu=\psi'(u)-a\Delta u,\qquad a>0,
\end{equation}
which results into the celebrated \emph{Cahn\pier{--}Hilliard system} \cite{CahnH1958JCP}. Yet, other choices have been considered: Novick-Cohen and Pego \cite{NovicP1991TAMS} and Plotnikov \cite{Plotn1994Passing} have investigated the \emph{viscous regularization}: 
\begin{equation}\label{eq:57}
\mu=\psi'(u)+\alpha\dot u,\qquad\alpha>0;
\end{equation}
furthermore, a combination of energetic and viscous regularization, namely, 
\begin{equation}\mu=\psi'(u)-a\Delta u+\alpha\dot u,
\nonumber
\end{equation} 
which leads to the so-called \emph{viscous Cahn\pier{--}Hilliard equation}, was derived by Novick-Cohen in \cite{Novic1988viscous} and analytically investigated by Elliott and Garcke in \cite{elliott1996cahn} and by Elliott and Stuart in \cite{EllioS1996JDE}, with computations being carried out in \cite{bai1995viscous}. \giu
Its vanishing-viscosity limit was studied in \cite{thanh2014passage}. More sophisticated generalizations of the Cahn--Hilliard system that still incorporate a viscous \pcol{contribution} have been proposed and investigated in \cite{miranville1999model, schimpaw}. \finegiu \elena \elefin{In the same spirit, a similar viscous regularization has been introduced in the paper \cite{BCDGSS}, which deals with phase separation in binary alloys driven by mechanical effects.} \finele
\medskip

\noindent\textbf{The elliptic regularization as a microforce balance.} Derivations of diffusion models \emph{\elena \`a \finele la} Cahn\pier{--}Hilliard  have been proposed by Gurtin in \cite{Gurti1996PD} and by Podio-Guidugli in \cite{Podio2006RM}. Central to these derivations are the following ingredients:
\begin{itemize} 
\item a system of microforces, distinct from Newtonian forces, which obey their own balance laws and whose power expenditure is associated to the evolution of the observable fields of interest --- in the present case, the concentration $u$; 
\item a collection of constitutive prescriptions, which relate microforces to the actual evolution of the system;
\item a dissipation principle, which sifts away thermodynamically inconsistent constitutive prescriptions.
\end{itemize}
Within this common framework, the elliptic regularization \eqref{eq:58} is a particular instance of the balance statement: 
\begin{equation}
\label{eq:17}
  \div{\boldsymbol\xi}+\pi+\gamma=0,
\end{equation}
an instance that arises on adopting the constitutive prescriptions:
\begin{equation}\label{eq:4}
  \boldsymbol\xi=a\nabla u,\qquad \pi=\mu-\psi(u),\qquad \gamma=0,
\end{equation}
for the vectorial \emph{microstress} $\boldsymbol\xi$ \pier{and} the two scalar-valued fields, $\pi$ and $\gamma$, respectively, the \emph{internal} and the \emph{external microforce}. 

Other constitituve prescriptions may be taken into consideration of course, provided that they are thermodynamically consistent. In this respect, two distinct options are offered in \cite{Gurti1996PD} and \cite{Podio2006RM}. In this paper we \pier{opt for} the former, where thermodynamical consistency is embodied by the dissipation inequality:
\begin{equation}\label{eq:3}
  \dot\phi\le (\mu-\pi)\dot u+{\boldsymbol\xi}\cdot\nabla\dot u-\mathbf h\cdot\nabla\mu,
\end{equation}
with $\phi$ the \emph{free-energy density}.

A standard argument \cite{ColemN1963ARMA} shows that consistency with \eqref{eq:3} rules out any constitutive dependence of free energy on the time derivative of $u$; accordingly, one assumes that the free energy and the concentration fields are related through a prescription of the form:
\begin{equation}\label{eq:1}
  \phi=\widehat\phi(u,\nabla u).
\end{equation}
A consequence of \eqref{eq:3} and \eqref{eq:1} is that, on introducing the \emph{equilibrium \pier{parts}}
\begin{equation}\label{eq:12}
  {\boldsymbol\xi}^{\rm eq}:=\frac{\partial\widehat\phi}{\partial\nabla u}(u,\nabla u)\qquad \textrm{and}\qquad \pi^{\rm eq}:=\mu-\frac{\partial\widehat\phi}{\partial u}(u,\nabla u),
\end{equation}
of, respectively, microstress and internal microforce, the test for consistency of a certain constitutive choice with the dissipation inequality \eqref{eq:3} boils down to verifying that, for whatever process, the \emph{non-equilibrium \pier{parts}}
\begin{equation}\label{eq:9}
  {\boldsymbol\xi}^{\rm ne}:={\boldsymbol\xi}-{\boldsymbol\xi}^{\rm eq},\qquad \pi^{\rm ne}:=\pi-\pi^{\rm eq}
\end{equation}
of microstress and internal microforce satisfy, together with the flux of diffusant $\mathbf h$, the \emph{reduced dissipation inequality}:
\begin{equation}\label{eq:13}
  0\le -\pi^{\rm ne}\dot u+{\boldsymbol\xi}^{\rm ne}\cdot\nabla\dot u-\mathbf h\cdot\nabla\mu.
\end{equation}
In particular, the constitutive prescriptions
 \eqref{eq:4} follow from \eqref{eq:12}--\eqref{eq:9} on taking
\begin{equation}
\widehat\phi(u,\nabla u)=\frac 12 a|\nabla u|^2+\psi(u),
\nonumber
\end{equation}
and on choosing
\begin{equation}\label{eq:8}
  \boldsymbol\xi^{\rm ne}=\mathbf 0,\qquad \pi^{\rm ne}=0,
\end{equation}
a choice consistent with  \eqref{eq:13}.
\medskip

\noindent\textbf{The viscous and the ``non smooth'' regularizations.} The viscous variant \eqref{eq:58} is arrived at in a similar fashion: first, we exclude microscopic contact interactions from the picture by letting $\boldsymbol\xi=\boldsymbol 0$, and we set to null the external microforce $\gamma$, so that the microforce balance \eqref{eq:17} reduces to
\begin{equation}\label{eq:64}
  \pi=0.
\end{equation}
\pier{Then}, consistent with this choice, we rule out the dependence of free energy on concentration gradient by letting
\begin{equation}\label{eq:67}
 \widehat\phi(u,\nabla u)=\psi(u),
\end{equation}
so that the second of \eqref{eq:12} specializes to
\begin{equation}\label{eq:66}
  \pi^{\rm eq}=\mu-\psi'(u);
\end{equation}
\elena As \pcol{to} the ``non-equilibrium part'', \elefin{the simplest constitutive choice is} \finele
\begin{equation}\label{eq:62}
  \pi^{\rm ne}=-\alpha\dot u,\qquad \alpha>0,
\end{equation}
so that, bearing in mind the second of \eqref{eq:9}, we recover \eqref{eq:58} from \eqref{eq:64} and \eqref{eq:66}.
\medskip

\elena Note that, actually, the above relation could be introduced in terms of a ``dissipation functional'' $\Phi$: the so-called pseudo-potential of dissipation by Moreau, which is a non-negative and convex functional, equal to zero for null dissipation. In particular, letting $\Phi(\dot u)=\frac \alpha 2 |\dot u|^2$, the non-equilibrium part could be introduced as
\begin{equation}\label{eq:62bis}
  \pi^{\rm ne}=-\frac{\partial\Phi}{\partial\dot u}.
\end{equation}
Note in particular that the assumption on $\Phi$ and \eqref{eq:62bis} (see \eqref{eq:62})  lead to 
\begin{equation}\label{eq:2}
\pi^{\rm ne}\dot u\le 0,
\end{equation}
which is important to ensure thermodynamical consistency. Indeed,  if \eqref{eq:2} and \eqref{eq:5} hold, then the reduced dissipation inequality \eqref{eq:13} is satisfied (recall that $\boldsymbol\xi^{\rm ne}=\boldsymbol 0$, which follows from $(\ref{eq:12})_1$, $(\ref{eq:9})_1$, and \eqref{eq:67}). 

Moreover, we could generalize \eqref{eq:62} by a suitable choice of $\Phi$. \finele Following a suggestion in \cite{Tomas2015?}, in this paper we replace \eqref{eq:62} with the following prescription:
\begin{equation}\label{eq:111}
  -(\pi^{\rm ne}+\alpha\dot u)\in\beta(\dot u),
\end{equation}
where, as before, $\alpha>0$ is a constant and $\beta:\mathbb R\rightrightarrows\mathbb R$ is a set-valued mapping whose graph is (maximal) monotone and containins the origin:
\begin{equation}\label{eq:10}
  0\in\beta(0),
\end{equation}
\elena as in the case $\beta=\partial \zeta$ is the subdifferential of a non-negative convex function $\zeta$, with $\zeta(0)=0$.\finele 

\elena As \pier{to} the assumption on $\beta$, \finele it is easily seen that the constitutive prescription \eqref{eq:111} is consistent with the dissipation inequality \eqref{eq:13}: indeed, owing to the monotonicity of $\beta$,
we have
\begin{equation}\label{eq:7}
  v_1\in\beta(w_1)\  \& \ v_2\in\beta(w_2)\quad\Rightarrow \quad(v_1-v_2)(w_1-w_2)\ge 0;
\end{equation}
thus, if the pair $(\pi^{\rm ne},\dot u)$ is compliant with \eqref{eq:111}, then we can take $v_1=-(\pi^{\rm ne}+\alpha\dot u)$ and $w_1=\dot u$ as tests in \eqref{eq:7}; meanwhile, \eqref{eq:10} entitles us to choose $w_2=v_2=0$; with these choices, \eqref{eq:7} yields $-(\pi^{\rm ne}+\alpha\dot u)\dot u\ge 0$, which entails  \elena \eqref{eq:2}.\finele

\medskip

\noindent\textbf{The system we investigate.} In order to assemble the system we study, it remains for us to: 

1) \pier{rewrite the balance equation \eqref{eq:11} in the light of} the constitutive prescription \eqref{eq:5} for the flux of diffusant, whence the partial differential equation:
\begin{equation}
  \dot u=m\Delta\mu.\label{eq:14}
\end{equation}

2) combine the microforce balance, in its form \eqref{eq:64}, with the second of \eqref{eq:9} and with the prescriptions \eqref{eq:66} and \eqref{eq:111}, whence the pointwise inclusion:
\begin{align}
\mu\in \psi'(u)+\alpha\dot u+\beta(\dot u),\label{eq:15}
\end{align}
where, with slight abuse of notation, we write $p+\beta(q)$ to denote the set $\{r\in\mathbb R:r-p\in \beta(q)\}$. 

We \elena finally complete \finele the system \pier{\eqref{eq:14}--\eqref{eq:15}}  with a prescription of the \emph{initial concentration} $u(0)$ (see \eqref{eq:18} below) and a time-dependent prescription of \pier{the} chemical potential on the boundary:
\begin{equation}\label{eq:31}
  \mu(\cdot,t)=\mu_\flat(\cdot,t)\quad\textrm{on }\, \pier{\Gamma := \partial \Omega}.
\end{equation}

Dirichlet-type boundary data involving chemical potential, rather than the flux of diffusant, are typical in applications. An example in the context of the mathematical modeling of hydrogen storage \cite{BonetCL2012NATMA,BCLkyoto, BonetFL2007AMO,roubicek2013thermomechanics} is the following: consider that the region $\Omega$ where diffusion takes place \elena represents \finele a body immersed in a gaseous reservoir at uniform pressure $p_g$ and temperature $T_g$, both of which may be possibly time dependent. The chemical potential of the diffusant in the reservoir is related to its pressure and temperature through the formula $\mu_g=\mu_0+RT_g\log(p_g/p_0)$ (see for instance \cite[Eq. 5.16]{AtkinP2006}), where $\mu_0$ is the \emph{standard chemical potential}, $R$ is the \emph{Avogadro constant}, and $p_0$ is the  \emph{standard pressure} (typically, $p_0=1{\rm bar}$). If local equilibrium prevails, $\mu$ is continuous across the boundary $\pier{\Gamma}$, and hence \eqref{eq:31} holds with $\mu_\flat(x,t)=\mu_g(t)$. Another example is provided by mechanical theories that describe the behaviour of a permeable elastic solid immersed in an incompressible fluid (see for instance \cite{DudaSF2010JMPS}); in these theories, the chemical potential of the fluid is given by $\mu_f=\mu_0+\nu(p_f-p_0)$, where $\nu$ is the molar volume of the solvent and $p_f$ is its pressure.
\medskip

\noindent\textbf{Energy and dissipation.} The precise assumptions we make on the coarse-grain free energy $\psi$ are stated in \eqref{eq:23} below. \elena Note in particular that we can allow  $\psi$ to be nonconvex  (a feature that, as already pointed out, allows for phase separation) and that we can include logarithmic type \pier{potentials}. On the other hand, we are not able to deal, e.g., with subdifferential of indicator functions of closed intervals.\finele

As to the mapping $\beta$, whose choice together with that of the constant $\alpha$ affects 
the dissipative structure of the system, we assume in \eqref{eq:17} below that it is the subdifferential of a \elena non-negative,\finele\  convex, lower semicontinous potential $\zeta$ 
(without setting any restriction on the growth of \pier{$\zeta$)} \elena with $\zeta(0)=0$
\finele. Besides the trivial \pier{case} $\beta=0$, which leads to the PDE considered in \cite{NovicP1991TAMS}, other possible choices are:
\begin{itemize}
\item [(i)] $\beta(r)=\beta_0\, \pier{\sign}(r)$, where $\beta_0>0$ and the sign graph is defined by
  \begin{equation}\label{eq:21}
\pier{\sign}(r)=\begin{cases}
\{+1\} &\textrm{ if }r>0,
\\
[-1,+1]&\textrm{ if }r=0,
\\
\{-1\} &\textrm{ if }r<0,
\end{cases}
\end{equation}
which, as we shall discuss below, \emph{may induce hysteresis}.
\item [(ii)] \elena the subdifferential of the indicator function $I_{[a,b]}$ of  a closed bounded interval  $[a,b]$ (in our computations we require $0\in[a,b]$) 
  \begin{equation}
    \beta(r)=\partial I_{[a,b]}(r)=
\begin{cases}
\{0\} &\textrm{ if }r\in (a,b),
\\
[0,+\infty)&\textrm{ if }r=b,
\\
(-\infty,0]&\textrm{ if }r=a,
\end{cases}
  \end{equation}
forcing \emph{the rate of change of concentration} to be bounded in the interval $[a,b]$.\finele 
\item [(iii)] the \elena subdifferential of the \finele indicator function of $[0,+\infty)$, namely,
 \begin{equation}
   \elena \beta(r)=\partial I_{[0,+\infty)}(r)=
\begin{cases}
\{0\} &\textrm{ if }r\in (0,+\infty),
\\
(-\infty,0]&\textrm{ if }r=0,
\end{cases}
\finele
  \end{equation}
which is a choice particularly interesting, for it entails that the concentration at a given point \emph{cannot decrease}, \pier{that is,  $\dot u$ in \eqref{eq:15} has to remain \elefin{non-negative}.} 
\end{itemize}
\medskip
%

\noindent\textbf{\pier{Remarks on} hysteresis.} 
For the viscous variant of the Cahn\pier{--}Hilliard system --- namely, the system that arises from \eqref{eq:57} --- it is known that hysteresis (in the sense of irreversibility \cite{Visin1994}) emerges in the \emph{vanishing\pier{-}viscosity limit}:
\begin{equation*}
\alpha\to 0,
\end{equation*}
provided that $\psi$ is a \emph{non-convex function} \cite{EvansP2004MMMAS,visintin2002forward}. It is not hard to construct examples showing that our constitutive assumptions lead to hysteresis as well, \emph{even if $\psi$ is convex}. In order to provide an illustration, we choose $\beta$ as in (i) above and we take $\psi(u)=\frac 12 ku^2$, with $k>0$. As a result, 
\pier{the system \eqref{eq:14}--\eqref{eq:15}} becomes:
\begin{subequations}\label{eq:63}
\begin{align}
&\pier{\partial_t u }=m\Delta\mu, \label{eq:63a}\\
&\mu\in \alpha\, \pier{\partial_t u }+\beta_0\,\pier{\sign}\left(\pier{\partial_t u }\right)+k u,\label{eq:63b}
\end{align}
\end{subequations}
\pier{where $ \pier{\partial_t}$ denotes the (partial) derivative with respect to $t$.} 
We supplement \eqref{eq:63a} with a boundary condition of the form
\begin{equation}\label{eq:61}
\mu(\cdot,t)=f\left(\dfrac t\tau\right)\quad\textrm{on }\pier{\Gamma},
\end{equation}
and we investigate the formal limit when the \emph{characteristic time} $\tau$ tends to infinity, which corresponds to the regime of a \emph{slowly-varying} chemical potential imposed at the boundary. 

On replacing $t$ with the dimensionless variable 
\begin{equation*}
s:=\frac t \tau,  
\end{equation*}
and on considering  that \pcol{the} $\beta$ \elefin{in (i)} is \emph{invariant under time reparametrization}, we can rewrite \eqref{eq:63}~\pier{as}
\begin{subequations}\label{eq:60}
\begin{align}
&\tau^{-1}\pier{\partial_s u }=m\Delta\mu,\label{eq:60a}\\
&\mu\in \tau^{-1}\alpha\,\pier{\partial_s u }+\beta_0\,{}\pier{\sign}\left(\pier{\partial_s u }\right)+k u.\label{eq:60b}
\end{align}
\end{subequations}
Formally, \pier{as} $\tau\to+\infty$, the parabolic system \eqref{eq:60} degenerates into the elliptic system:
\begin{subequations}\label{eq:59}
\begin{align}
&0=m\Delta\mu,\label{eq:59a}\\
&\mu\in \beta_0\,\pier{\sign}\left(\pier{\partial_s u }\right)+k u.\label{eq:59b}
\end{align}
\end{subequations}
Having stipulated with \eqref{eq:61} that  $\mu$ is spatially \elena constant \finele  on the boundary, the homogeneous elliptic equation \eqref{eq:59a} entails that $\mu$ is spatially uniform in the bulk:
\begin{equation}
  \mu(\cdot,s)=f(s)\quad\textrm{in }\Omega.
\end{equation}
 As a consequence, the concentration field satisfies the following differential inclusion:
\begin{equation}\label{eq:54}
  f(s)\in \beta_0\,\pier{\sign}\left(\pier{\partial_s u }\right)+k u,
\end{equation}
which is known to exhibit hysteresis (actually, it reproduces the well-known \emph{stop operator}, see e.g. \cite{brokate1996hysteresis, Visin1994}).
\medskip

\noindent\elefin{{\bf The initial boundary value problem.}} We consider the evolution in a smooth domain $\Omega$ in the expanse of time $(0,T)$. We suppose that at time $t=0$ the concentration field be given by a prescribed function $u_0(x)$, $x\in\Omega$. We also suppose that a time-dependent chemical potential $\mu_\flat(x,t)$ be prescribed on for all $x\in\pier{\Gamma}$ at all times $t\in[0,T]$. 

At each particular time, we harmonically extend $\mu_\flat$ to the interior of $\Omega$ \elena(still \pcol{denoting} by $\mu_\flat$ the harmonic extension)\finele\ and we introduce the characteristic time and lengthscale $T_0=\alpha$ and $L_0=\sqrt{\frac{m}{\alpha}}$,\footnote{We assume that all energies densities per unit volume are renormalized to a reference value, and so are dimensionless.} \pier{as well as} the new variables and functions:
\begin{equation}
  \widetilde t=\frac t {T_0},\qquad \widetilde x=\frac x {L_0},\qquad \widetilde \mu=\mu-\mu_\flat,\qquad \widetilde\beta(r)=\beta(r/T_0).
\end{equation}
We express the system in \pier{terms of these} new variables and \emph{we henceforth drop tildas}, so as to obtain the following problem:
\begin{subequations}\label{eq:16}
\begin{equation}
\begin{aligned}
\left.
\begin{array}{l}
  \pier{\partial_t u }=\Delta \mu,\\
  \mu=\pier{\partial_t u }+\xi+\mu_\flat+\psi'(u)\\
  \xi\in\beta(\pier{\partial_t u })
\end{array}
\right\}\quad\textrm{in }\, \pier{\Omega \times} (0,T)
\end{aligned}
\end{equation}
 with the boundary condition
\begin{equation}\label{eq:19}
  \mu=0\quad\textrm{ \pier{on}  }\,\pier{\Gamma \times } (0,T),
\end{equation}
and \elefin{the} initial condition
\begin{align}\label{eq:18}
  u(\cdot,0)=u_0\quad\textrm{in } \Omega.
\end{align}
\end{subequations}

\elena
Before proceeding, let us \pier{set a precise outline} of the paper. In the next section, we introduce notation, assumptions on the data of the problem and state the main existence and uniqueness result, which is complemented by \pier{the}  continuous dependence of \pier{the} solution with respect to the data of the problem. In Section~\ref{sec:basic-estimates} we  proceed by exploiting the a priori estimates on the solutions of the system we use to prove existence and regularity. In Section~\ref{sec:cont-depend-data-1} we \pier{show the} continuous dependence estimates on solutions. 
\elefin{Finally, in Section~\ref{sec:5pier}, we provide a detailed proof of the existence and the uniqueness of the solution. To this aim, we first establish our result for a ``regularized'' version of the system, obtained by a suitable truncation of the coarse-grain free energy mapping; then, owing to the estimates carried out in Section~\ref{sec:basic-estimates}, and by making use of a maximum--principle argument, we establish our result for the original problem.} 
\finele
\section{Notation, assumptions, and main results}
\setcounter{equation}{0}
\elena
Before stating the problem we are dealing with and the main existence result, let us make precise our assumptions on the data and the notation we use.
In the sequel, $\Omega$ is a bounded \pier{smooth domain in~$\mathbb{R}^3$ with smooth
boundary $\Gamma$}. \pier{We introduce the spaces}
$$
  H := L^2(\Omega), \qquad V := H^1_0(\Omega), \qquad W := H^2(\Omega)\cap H^1_0(\Omega).
$$
We \pier{endow} $H$, $V$, and $W$
with their usual scalar products and norms, and use a \pcol{self-explanatory}
notation, like $\|\cdot\|_V$.  For the sake of simplicity, the same symbol will be
used both for a space and for any power of it.  We note that the norms
$\|v\|_V$ and $\|\nabla v\|_H$ are equivalent for $v\in V$, thanks to
the Poincar\'e inequality. In addition, let us point out that, after identifying $H$ with its dual, the triplet $(V,H,V')$ is a Hilbert triplet (where $V'$ coincides with the Sobolev space $H^{-1}(\Omega)$). 
Hence, we use the notation $\langle\cdot,\cdot\rangle$ for the duality pairing between $V'$ and $V$. \pier{Given a final time~$T>0$, we~set}
$$
  Q := \Omega \times (0,T), \quad \pier{\Sigma := \Gamma \times (0,T).}
$$

As far as the data of the problem are concerned, we \pier{assume} that the domain of $\psi$ is an open interval  $(a,b)\subseteq\mathbb{R}$, \pcol{where $a$ and $b$ could be taken equal to $-\infty$ and $+\infty$, respectively}%
\pier{, and require that}\finele
\begin{subequations}\label{eq:23}
\begin{gather}
  \elena\psi\in C^2(a,b),\label{eq:23a}\finele\\
 \gius{\psi(r)\ge 0} \quad \pier{\hbox{for all } \, r\in (a,b)} ,\label{eq:23b}\\ 
 \pier{ \lim_{r\to a^+}\psi'(r)=-\infty , \quad \lim_{r\to b^-}\psi'(r)=+\infty,} \label{eq:23d}\\
 \psi''(r)\ge \gius{-K_1}  \quad \pier{\hbox{for all } \, r\in (a,b)} , \label{eq:23c}
\end{gather}
for some positive \elena constant $K_1$;\finele\ \pcol{note that in \eqref{eq:23d}} $a^+ $ has to become $-\infty$ if $a= -\infty$ and $b^-$ reduces to $+\infty$ if $b=+\infty$.
\elena 
\pier{For the initial concentration $u_0$ we suppose that
\begin{align}\label{eq:6}
u_0\in H, \quad \elefin{\exists} \, a_0>a, \ \, b_0 < b\,  \hbox{ such that }\,  a_0\leq u_0 (x) \leq b_0 \, \hbox{ for a.a. } x\in\Omega,
\end{align}
whence both $ u_0 $ and $\psi'(u_0)$ lie in $ L^\infty(\Omega)$.
Concerning the \elefin{known datum}  $\mu_\flat $, we assume that}
\begin{equation}\label{eq:24}
  \mu_\flat\in H^1(0,T;\pier{H})\cap \pier{L^\infty(Q)}.
\end{equation}
Finally, \finele as to the \pier{(possibly)} multi-valued mapping $\beta$, we \pier{let} \elefin{(see \cite{barbu})}
\begin{gather}
\hbox{$\beta=\partial\zeta$\pier{, with $\zeta:\mathbb R\to[0,+\infty]$ \emph{convex} and}}\nonumber \\ 
\hbox{\pier{\emph{lower-semicontinuous}, such that $\zeta(0)=0$.}\label{eq:17ele}}
\end{gather}
\end{subequations}

\pier{Let us specify a} \elena weak formulation of the problem in the set of the Hilbert triplet $(V,H,V')$. \pier{Let us define} the operator $A:V\rightarrow V'$, corresponding 
to the ``weak realization'' of the Laplace operator  $-\Delta$ (combined with homogeneous Dirichlet boundary condition) in the duality  \elefin{between $V'$ and $V$},  by letting \finele
\begin{align}\label{eq:20}
{}\langle Av_1,v_2\rangle:=\int_\Omega \nabla v_1\cdot\nabla v_2 \,{\rm d}x, \quad \pier{v_1,v_2\in V,}
\end{align}
\elena
and \pier{specify} its inverse $A^{-1}:V'\rightarrow V$, such that for $w,z\in V'$ \pier{and} $v\in V$
\begin{equation}\label{eq:20bis}
\langle Av,A^{-1}w\rangle=\langle w,v\rangle,\quad \langle w,A^{-1}z\rangle=\langle z,A^{-1}w\rangle=\int_\Omega\nabla(A^{-1}w)\cdot\nabla(A^{-1}z).
\end{equation}
\pier{We introduce a norm in $V'$, denoted by $\|\cdot\|_*$, which is equivalent to the usual \pier{one,}}
\begin{equation}\label{eq:20tris}
\|w\|^2_*=\langle w,A^{-1}w\rangle, \pier{\quad w\in V'.}
\end{equation}

\noindent{\bf Definition of solutions}. We  say that a triplet $(u,\xi,\mu)$ is a weak solution to \pier{the problem \eqref{eq:16} if}
\begin{subequations}\label{pier2} 
\begin{gather}
  u\in C^1([0,T];H),\quad
  \xi\in C^0([0,T];H),\quad
  \mu\in C^0([0,T];\pier{V}),\label{pier2a}\\
  \pier{\psi'(u)\in H^1(0,T;H)} \label{pier2b} 
\end{gather}
\end{subequations}
and  the following equations are satisfied:
\begin{subequations}\label{eq:22}
\begin{gather}
\pier{\partial_t u }(t)+A\mu(t)=0\quad \textrm{in }V',\, \hbox{ for all }t\in[0,T],\label{eq:6-1}\\
\pcol{\mu(t) = (\pier{\partial_t u }+\xi+\mu_\flat+\psi'(u))(t) \quad\textrm{ a.e. in }\Omega,\, \hbox{ for all }t\in[0,T],}\label{eq:6-2}\\
\pcol{\xi(t)\in\beta(\pier{\partial_t u }(t) )\quad\textrm{ a.e. in }\Omega, \, 
\hbox{ for all }t\in[0,T],}\label{eq:6-3}\\
u(0)=u_0\quad\textrm{a.e. in } \Omega.\label{eq:6-4}
\end{gather}
\end{subequations}\finele
\begin{theorem}[Existence and uniqueness]\label{Th1}
Under \pier{the} assumptions \eqref{eq:23}, there exists a \elena unique  solution to 
\pier{the problem~\eqref{eq:16}}, in the sense we have specified above. In addition, it results that
\begin{equation}
u\in L^\infty(Q)
\label{pier10}
\end{equation}
and \pier{the solution is strong, that is, $\mu \in C^0([0,T];W)$ and equation \eqref{eq:6-1} can be replaced~by 
\begin{equation}
\label {pier1}
\pcol{ \pier{\partial_t u }(t) - \Delta \mu (t) =0 \quad \hbox{a.e. in } \Omega , \,
\hbox{ for all }t\in[0,T].}
\end{equation} 
}\finele Moreover, \pier{a continuous dependence on the data holds: namely,} if $(u_1,\xi_1,\mu_1)$, $(u_2,\xi_2,\mu_2)$ are two solution triplets corresponding to the initial data $u_{01}$, $u_{02}$ and bulk data $\mu_{\flat1}$,  $\mu_{\flat2}$, \elena respectively, \finele then their difference satisfies
\elena
\begin{align}\label{contdep}
 & \|u_1-u_2\|_{C^1([0,T];H)}+\|\xi_1-\xi_2\|_{C^0([0,T];H)}+\|\mu_1-\mu_2\|_{C^0([0,T];W)}\nonumber\\
 &\le \pcol{R} \,\bigl(\|\mu_{\flat1}-\mu_{\flat2}\|_{C^0([0,T];H)}+\|u_{01}-u_{02}\|_{H}\bigr)
\end{align}
\finele
\pier{for some constant \pcol{$R$} depending only on the structural assumptions stated in \eqref{eq:23}.}
\end{theorem}
\section{Basic estimates}
\label{sec:basic-estimates}
\setcounter{equation}{0}
\elena In this section, for reader's convenience, before proving Theorem \ref{Th1}, we \pier{recover} the a priori estimates we can derive on the solutions of \finele  \eqref{eq:22}. \elena Note that we are using the same notation $C$ for possible different positive constants depending only on the data of the problem.\finele
\medskip

\noindent \textbf{\elena Energy \finele estimate.} We \elena first show that system \eqref{eq:22} admits in a natural way a so-called ``energy estimate''. Indeed, once $u$ and $\mu$ \pier{are solution components for the problem~\eqref{eq:16},} we are allowed to \finele test \eqref{eq:6-1} by $\mu$ and \eqref{eq:6-2} by $\pier{\partial_t u }$. \pier{Then, we combine the resulting equations and integrate by parts in time over $(0,t)$; by exploiting H\"older's and Young's inequalities and using the chain rule and the smoothness of $\psi$, we \elefin{find} that}
\elena
\begin{align}\label{eq:46}
 & \int_\Omega \psi(u(t))\,{\rm d}x+\int_0^t\!\! \int_\Omega (|\pier{\partial_t u }|^2+\xi \, \pier{\partial_t u }+|\nabla\mu|^2)
 \,{\rm d}x\pier{{\rm d}s}\nonumber\\
 &\le \int_\Omega \psi(u_0)\,{\rm d}x+\int_0^t\!\!\int_\Omega |\mu_\flat \, \pier{\partial_t u }|{\rm d}x\pier{{\rm d}s},\nonumber\\
 &\le \int_\Omega \psi(u_0)\,{\rm d}x+\frac 12\int_0^t\!\!\int_\Omega |\mu_\flat|^2{\rm d}x\pier{{\rm d}s}+\frac 12\int_0^t\!\!\int_\Omega |\pier{\partial_t u }|^2{\rm d}x\pier{{\rm d}s}.
\end{align}
\pier{Owing to the} monotonicity of $\beta$ and the fact that $0\in\beta(0)$, we \pier{have that}
\begin{equation}
\int_0^t\!\! \int_\Omega \xi \, \pier{\partial_t u }\,{\rm d}x\pier{{\rm d}s}\geq0.
\end{equation}
\elefin{As \eqref{pier2a} holds}, we also point out that\  $u(t) = u_0 + \int_0^t\partial_t {u}(s){\rm d}s,$ whence
$$ \| u(t)\|_H^2\le 2 \| u_0\|^2_{H}+ 2T \!\int_0^t\|\pier{\partial_t {u}}(s)\|_H^2\,\pier{{\rm d}s}, $$ 
thanks to the H\"older inequality. Thus, \finele by virtue of the Poincar\'e inequality \pier{\gius{and the \pier{nonnegativity} of $\psi$ as well,} \eqref{eq:46} and \eqref{eq:23} yield the estimate}
\begin{equation}\label{stimaenergia}
 \|\psi(u)\|_{L^\infty(0,T;L^1(\Omega))}+ \|u\|_{H^1(0,T;H)}+\|\mu\|_{L^2(0,T;V)}\le C.
\end{equation}
\elena Note that \pier{\eqref{stimaenergia} entails in particular that $u$ lies between $a$ and  $b$} almost everywhere in $Q$.%
\finele
\pier{Hence, recalling that $\Omega$ and $\Gamma$ are smooth enough and that \elefin{$\mu$ has null trace on $\Gamma$}, by \elena a comparison in \finele (\ref{eq:6-1}) and standard 
elliptic regularity estimates} \elefin{we obtain, in addition,} \pier{\eqref{pier1} and}
\begin{equation}
  \|\mu\|_{L^2 (0,T;W)}\le C.
\end{equation}
\medskip

\noindent\textbf{$\mathbf{L^\infty}$ estimate for chemical potential.} \pier{Let us first introduce $\mu_0\in V$ as the unique solution of the nonlinear elliptic problem
\begin{equation}\label{servemuzero}
  A\mu_0+(I+\beta)^{-1}(\mu_0-\mu_\flat(0)-\psi'(u_0))=0.
\end{equation}
Indeed, note that  $\mu_\flat(0)+ \psi'(u_0)\in H$ by \pcol{\eqref{eq:24}, \eqref{eq:6}} \elefin{and \eqref{eq:23a}}: then, \elefin{concerning} the sum of the \giu two \finegiu maximal monotone operators $A$ (when restricted to $W$ with values in $H$)~and  
$$v\mapsto (I+\beta)^{-1}(v-\mu_\flat(0)-\psi'(u_0))$$ 
we can apply \cite[Cor.~1.3, p.~48]{barbu}, which ensures that the sum is maximal monotone and surjective thanks to the Lipschitz continuity of the second operator and the coerciveness of~$A$.
Moreover, the uniqueness of $\mu_0\in W$ solving 
\eqref{servemuzero} follows from the strong monotonicity of~$A$. Clearly, from \eqref{pier2} and \eqref{eq:22} we have that $\mu(0)= \mu_0 $ and
\begin{equation}\label{eq:41}
\pier{\partial_t u }(0)=\Delta\mu_0=:\pier{u'_0}\in H\quad\textrm{and}\quad \xi(0)=\mu_0- \pier{u'_0}-\mu_\flat(0)-\psi'(u_0)=:\xi_0
\in H.
\end{equation}
We emphasize that $\xi_0$ satisfies $\xi_0 \in \beta(\pier{u'_0})$ almost everywhere in $\Omega$. Next, in order to show the $L^\infty$ estimate for $\mu$, let us fix $t\in (0,T]$ and, for $n\in\NN$, set $\tau_n= t/n$, \elefin{$t_i^n=i\tau_n$, $i=0, 1, \ldots, n$. For typographical convenience, we henceforth omit the dependence of $t_i^n$ on $n$.} In view of \eqref{pier2}, 
equations \eqref{eq:6-1} and \eqref{eq:6-2} hold at the times $t_i$:
\begin{gather}
\pier{\partial_t u }(t_i)+A\mu(t_i)=0\quad \textrm{in }V',\label{pier3}\\
\mu(t_i)= \pier{\partial_t u }(t_i)+\xi(t_i)+(\mu_\flat+\psi'(u))(t_i)\quad\textrm{in }H,\label{pier4}
\end{gather}
for $i= 1, \ldots, n$. Test \eqref{pier3} by $ \mu (t_i) - \mu (t_{i-1})$ and the difference of  equalities \eqref{pier4} at the steps $i$ and $i-1$ by $ \partial_t u (t_i)$. Then, by 
combining the results is not difficult to check that 
\begin{align}
&\nonumber \frac 12 \int_\Omega |\nabla\mu(t_i )|^2{\rm d}x - \frac 12 \int_\Omega |\nabla\mu(t_{i-1} )|^2{\rm d}x + \frac 12 \int_\Omega |\nabla\mu(t_i ) - \nabla \mu(t_{i-1} )|^2{\rm d}x\\
&\nonumber {}+ \frac 12 \Vert \partial_t u (t_i) \Vert_H^2   - \frac 12 \Vert \partial_t u (t_{i-1}) \Vert_H^2
+ \frac 12 \Vert  \partial_t u (t_i)- \partial_t u (t_{i-1}) \Vert_H^2\\
\label{pier5} &{}+ (\xi(t_i) - \xi(t_{i-1}), \partial_t u (t_i))+ ((\mu_\flat+\psi'(u))(t_i) - (\mu_\flat+\psi'(u))(t_{i-1}), \partial_t u (t_i)) =0 ,
\end{align}
where $(\cdot, \cdot) $ denotes the scalar product in $H$. Now, by the 
properties of the subdifferential $\beta=\partial\zeta$ (which is a maximal monotone graph) 
it turns out that the inclusion (see~\eqref{eq:6-3} and \eqref{eq:17ele}) 
$$\xi(t_i) \in \partial \zeta( \partial_t u (t_i)) $$
can be rewritten as
\begin{equation}\label{pier6}
  \pier{\partial_t u (t_i) }\in \partial\zeta^*(\xi(t_i))
\end{equation}
almost everywhere in $\Omega$,
where $\zeta^*(w):=\sup_{v\in\mathbb R} (v\, w -\zeta(v))$, $w\in\RR$,  is the 
Legendre-Fenchel transform of $\zeta$,  and its subdifferential $\partial\zeta^*$ coincides with $\beta^{-1}$, the 
inverse graph of $\beta$. Then, using the definition of subdifferential, it is straightforward to 
infer that 
$$  
(\xi(t_i) - \xi(t_{i-1}), \partial_t u (t_i)) \geq \int_\Omega\zeta^*(\xi(t_i)){\rm d}x  - \int_\Omega\zeta^*(\xi(t_{i-1})){\rm d}x .
$$
Hence, \elefin{on performing summation in \eqref{pier5} for $i=1, \ldots, n$}, we plainly deduce that
\begin{align}
&\nonumber \frac 12 \int_\Omega |\nabla\mu(t)|^2{\rm d}x 
+ \frac 12 \Vert \partial_t u (t) \Vert_H^2 + \int_\Omega\zeta^*(\xi(t)){\rm d}x \\
&{}+ \sum_{i=1}^n \tau_n \bigl( \frac{(\mu_\flat+\psi'(u))(t_i) - 
(\mu_\flat+\psi'(u))(t_{i-1})}{\tau_n}, \partial_t u (t_i)\bigr)\nonumber\\
&{}\leq \frac 12 \int_\Omega |\nabla\mu_0|^2 {\rm d}x + \frac 12 \Vert u'_0 \Vert_H^2 + \int_\Omega\zeta^*(\xi_0){\rm d}x .  \label{pier7} 
\end{align}
Owing to \eqref{eq:24} and \eqref{pier2b}, it is a standard matter to infer that
$$\sum_{i=1}^n \tau_n \bigl( \frac{(\mu_\flat+\psi'(u))(t_i) - 
(\mu_\flat+\psi'(u))(t_{i-1})}{\tau_n}, \partial_t u (t_i)\bigr)\to \int_0^t\!\!\int_\Omega \partial_t (\mu_\flat+\psi'(u)) \pier{\partial_t u }\,{\rm d}x\pier{{\rm d}s }$$
as $n\to \infty $, and assumptions \eqref{eq:23c} and \eqref{eq:24} easily yield 
\begin{equation*}
\int_0^t\!\!\int_\Omega \partial_t (\mu_\flat+\psi'(u)) \pier{\partial_t u }\,{\rm d}x\pier{{\rm d}s } \geq - \|\pier{\partial_t \mu_\flat }\|_{L^2(0,T;H)} \|\pier{\partial_t u }\|_{L^2(0,T;H)}
 -K_1\|\pier{\partial_t u }\|^2_{L^2(0,T;H)} \geq - C,
\end{equation*}
the last inequality being due to the previous estimate \eqref{stimaenergia}. As $\xi_0\in \partial \zeta (u'_0)$ and consequently $\zeta^* (\xi_0) + \zeta (u'_0)= \xi_0\, u'_0$ almost everywhere in $\Omega$, from \eqref{eq:17ele} it follows that
\begin{equation}\label{eq:30}
 \int_\Omega \zeta^*(\xi_0){\rm d}x\le \|\xi_0\|_{H}\|\pier{u'_0}\|_{H}-\int_\Omega \zeta(\pier{u'_0}){\rm d}x\le C.
\end{equation}
On the other end, it is easy to check that $\zeta^*$ is non-negative, whence $$\int_\Omega\zeta^*(\xi(t)){\rm d}x\geq0. $$
Then, passing to the limit as $n\to \infty$ in \eqref{pier7} and exploiting the
previous remarks, we find out that 
\begin{equation}\label{eq:44}
 \|\mu(t)\|^2_{V}+ \|\pier{\partial_t u }(t)\|^2_{H}\le C \quad \hbox{for all }\, t \in [0,T].
\end{equation}
At this point, by comparison in \eqref{eq:6-1} and thanks to well-known elliptic regularity 
results combined with the Sobolev embedding $W\subset L^\infty(\Omega)$, we recover
\begin{equation}
\|\mu\|_{L^\infty(0,T;W)}+\|\mu\|_{L^\infty(Q)}\le C.
\end{equation}
In} particular, because of the $L^\infty$ boundedness of $\mu_\flat$ \pier{postulated} in Assumption \eqref{eq:24}, we infer that \elena
there exists a constant $M$ such that \finele
\begin{equation}\label{eq:49}
  \|\mu-\mu_\flat\|_{L^\infty(Q)}\le M .
\end{equation}
\medskip

\noindent\textbf{$\mathbf{L^\infty}$ estimate for concentration.} We combine \pier{\eqref{eq:6-2}--\eqref{eq:6-4}} to obtain the Cauchy problem in $H$:
\begin{subequations}\label{eq:32}
\begin{gather}
  \pier{\partial_t u }(t)=(I+\beta)^{-1}(\mu(t)-\mu_\flat(t)-\psi'(u(t))) \label{pier8}
  \pier{,\quad t\in [0,T],}\\
  u(0)=u_0. \label{pier9}
\end{gather}
\end{subequations}
\gius{Now, by \pier{assumptions \eqref{eq:23d}} and \eqref{eq:6} there exists \pier{two} constants \elena  $k_*,k^*\in (a,b)$ \finele  such that}
\begin{subequations}\label{eq:33}
\begin{gather}
  \gius{\psi'(r)\ge M}\quad\textrm{for all }r\ge k^*,\label{eq:33b}\\ 
  \gius{\psi'(r)\le -M}\quad\textrm{for all }r\gius{\le k_*,}\label{eq:33c}\\
  \pier{k_*\le a_0\le u_0 (x)\le b_0\le  k^*\quad\hbox{for a.a. } x\in\Omega.}\label{eq:33a}
\end{gather}
\end{subequations}
\elena We test  \eqref{pier8} by $(u-k^*)^+ $ and $ - (u-k_*)^-$, then integrate over $(0,t)$. \finele 
\giu We n\finegiu ote that $(I+\beta)^{-1}(r)$ has the same sign of $r$ \pier{and the right hand side of \eqref{pier8} is nonpositive if $u\ge k^*$ and nonnegative if $u\le k_*$, thanks to \eqref{eq:33} and \eqref{eq:49}. Then,} \elena after \pier{integration} by parts in time, it is a standard matter to infer that\finele
\begin{equation}\label{eq:28}
  k_*\le u\le k^*\quad \textrm{a.e. in }Q,
\end{equation}
which entails \eqref{pier10}.
\medskip
\section{Continuous dependence on the data}
\label{sec:cont-depend-data-1}
\setcounter{equation}{0}
%
%
%
Consider a pair of data  \pier{$\{  u_{0i}, \, \mu_{\flat i}\} $, $i=1,2$, fulfilling \pier{\eqref{eq:6}, \eqref{eq:24}} and let \elena $(u_i,\xi_i,\mu_i)$,
$i=1,2$, 
\finele  be the corresponding solutions. Then, the triplet $(\bar u,\bar\xi,\bar\mu)$, with $\bar u:=u_1-u_2$, $\bar\xi:=\xi_1-\xi_2$, $\bar\mu:=\mu_1-\mu_2$,  satisfies (cf.~\eqref{eq:22} and \eqref{pier1})
\begin{subequations}\label{eq:34}
\begin{gather}
\pcol{\pier{\partial_t {\bar u }}(t) - \Delta \bar \mu (t) =0
\quad \hbox{a.e. in } \Omega , \, \hbox{ for all }t\in[0,T],}
\label{eq:34a}\\
\pcol{\bar\mu (t)= (\pier{\partial_t \bar{ u }}+\bar\xi+\bar\mu_\flat+\psi'(u_1)-\psi'(u_2))(t)\quad \hbox{a.e. in } \Omega , \, \hbox{ for all }t\in[0,T],}\label{eq:34b}\\
\pcol{\xi_i(t)\in\beta(\partial_t u_{i}(t) )\quad \hbox{a.e. in } \Omega , \, \hbox{ for all }t\in[0,T],} \, \  i=1,2, \label{eq:34c}\\
\bar u(0)=\bar u_ 0\quad\textrm{a.e. in } \Omega,\label{eq:34d}
\end{gather}
\end{subequations}
where $\bar u_0=u_{01}-u_{02}$ and $\bar\mu_\flat=\bar\mu_{\flat 1}-\bar\mu_{\flat 2}$.}

\pier{In view of the regularities in \eqref{pier2a}, we can test \eqref{eq:34a} by $\bar\mu \pcol{(t)}$, \eqref{eq:34b} by $ \pier{\partial_t \bar{u }}\pcol{(t)}$ and add the resulting equations. In particular, by virtue of \eqref{eq:34c} and by the monotonicity of $\beta$, we have that  $\int_\Omega \bar\xi(t)\pier{\partial_t {\bar u}}(t)\ge 0$ and consequently}
\begin{equation}
  \int_\Omega|\nabla \bar \mu(t)|^2\,{\rm d}x+\int_\Omega |\pier{\partial_t {\bar u}}(t)|^2\,{\rm d}x\le \int_\Omega (|\bar\mu_\flat(t)|+|\psi'(u_1(t))-\psi'(u_2(t))|)|\pier{\partial_t {\bar u}}(t)|\,{\rm d}x.
\end{equation}

\pier{Now, since $\psi$ is twice continuously differentiable, its derivative is locally Lipschitz-continuous. Moreover, by the estimate \eqref{eq:28}, it turns out that both solutions $u_i$ stay in a bounded interval $J$. Consequently, we have that $|\psi'(u_1)-\psi'(u_2)|\le \Vert \psi''\Vert_{L^\infty (J)} |\bar u|$ and, by Young's and Poincar\'e's  inequalities, we infer 
\begin{equation}\label{eq:35}
 \elena\|\bar\mu(t)\|^2_V\finele+ \|\pier{\partial_t {\bar u}}(t)\|_{H}^2\le C\bigl( \|\bar\mu_\flat(t)\|^2_H+\|\bar 
 u(t)\|^2_{H} \bigr).
\end{equation}
Given that $\|\bar u(t)\|_H^2\le C\left(\|\bar u_0\|^2_{H}+\int_0^t\|\pier{\partial_t {\bar u}}(s)\|_H^2\,\pier{{\rm d}s}\right)$, from \eqref{eq:35} it follows that
\begin{align*}
&\|\bar \mu(t)\|^2_V+  \|\pier{\partial_t {\bar u}}(t)\|_{H}^2  \nonumber \\
&\le C\left(\|\bar\mu_\flat(t)\|_H^2+\|\bar u_0(t)\|^2_{H}+\int_0^t\|\pier{\partial_t {\bar u}}(s)\|_H^2\,\pier{{\rm d}s}\right)\quad\textrm{for all }t\in[0,T].
\end{align*}
Thus, an application of the Gronwall-Bellmann inequality yields
\elena\begin{equation}\label{eq:36}
\|\bar\mu\|_{C^0([0,T];V)}+ \|\pier{\partial_t {\bar u}}\|_{C^0([0,T];H)}\le C\left(\|\bar\mu_\flat\|_{C^0([0,T];H)}+\|\bar u_0\|_{H}\right).
 \end{equation}
Then, the analogous estimates for $A\bar\mu$ \finele (and consequently for $\bar\mu$ in $C^0 ([0,T];W)$) and $\bar\xi$ in $C^0([0,T];H)$ follow from \eqref{eq:36} by a comparison in \eqref{eq:34a} and \eqref{eq:34b}, which helps us to conclude the proof of \eqref{contdep}.} \pier{Of course, \eqref{contdep} implies in particular the uniqueness of the solution to the problem \eqref{eq:22}.}
\section{Existence of solutions}
\label{sec:5pier}
\setcounter{equation}{0}

\elena
In this \pier{section},
we give details on the proof of the existence \pier{of the} solution to our problem. We use a contracting argument. First let us make a truncation of the function $\psi$, which allows us to exploit the above a priori bounds on the solutions to \pier{\eqref{eq:22}}. \finele \gius{Let $k_*$ and $k^*$ be \pier{two constants fulfilling~\eqref{eq:33}}. It is not hard to check that, thanks to the assumption \eqref{eq:23d}, there exist constants $K_*$ and $K^*$ such that}
\begin{equation}
  (k_*,k^*)\subseteq (K_*,K^*)
\end{equation}
and
\begin{equation}
  \psi''(K_*)\ge 0,\qquad \psi''(K^*)\ge 0.
\end{equation}
We introduce the following truncation of $\psi$:
\begin{equation}  
 \psi_*(r)=\begin{cases}
\psi(r)\quad &\textrm{if }\, K_*\le r\le K^*,\\
\psi(K^*)+\psi'(K^*)(r-K^*)+\frac 12 \psi''(K^*)(r-K^*)^2\ & \textrm{if }\, r> K^*,\\
\psi(K_*)+\psi'(K_*)(r-K_*)+\frac 12 \psi''(K_*)(r-K_*)^2\ & \textrm{if }\, r< K_*,\\
\end{cases}
\end{equation}
and we denote by
\begin{equation}\label{eq:52}
  L:=\max_{r\in\mathbb R}|\psi''_*(r)|<+\infty
\end{equation}
the Lipschitz constant of its derivative $\psi_*'$. We note on passing that the truncated function $\psi_*$ \pier{satisfies the assumptions (\ref{eq:23a})--(\ref{eq:23c})
with $(a,b)= (-\infty, +\infty)$}. In particular,  \pier{the bound from below\eqref{eq:23c} holds for $\psi''_*$ with the same constant $- K_1$ as for $\psi''$}.

Next, we consider the set
\begin{equation}
  S:=\{v\in \elena C^0\finele([0,T];H):\, v(0)=u_0\}
\end{equation}
and we introduce the map $\mathcal F:S\to S$ which to every $v\in S$ associates $u=\mathcal F(v)$ defined by
\begin{equation}\label{eq:42}
 u(t)=\mathcal F(v)(t):=u_0+\int_0^t (I+\beta)^{-1}\left(\mu(s)-\mu_\flat\pier{(s)} -\psi_*'(v(s))\right)\pier{{\rm d}s},\quad t\in[0,T],
\end{equation}
where $\mu(t)$ denotes the unique element of $V$ that solves the \pier{nonlinear elliptic equation}
\begin{equation}\label{punto2}
  A\mu(t)+(I+\beta)^{-1}(\mu(t)-\mu_\flat(t)-\psi_*'(v(t))=0,\quad t\in[0,T].
\end{equation}

\pier{Before proceeding, let us comment on the existence of a unique $\mu \in  C^0([0,T];V)$ satisfying \eqref{punto2} for some $v$ fixed in $S$. First, we recall \eqref{eq:24}) and observe that $\psi_*'$ is Lipschitz continuous, so that  the function $t \mapsto \mu_\flat(t)-\psi_*'(v(t))$ is continuous from $[0,T]$ to $H$. Then, for all $t\in [0,T]$ there exists a
unique $\mu(t) $ fulfilling \eqref{punto2}: this can be shown arguing as for   \eqref{servemuzero} and using \cite[Cor.~1.3, p.~48]{barbu}. Moreover, as 
$(I+\beta)^{-1}$ is monotone and Lipschitz continuous, it is not difficult to check that 
$\mu \in  C^0([0,T];V)$. Once $\mu$ is found, the function $u=\mathcal F(v)\in S $ is completely determined from  \eqref{eq:42}.}

\elena Eventually, our aim consists in applying a fixed point argument\pier{: indeed, we will see that} any fixed point for the operator $\mathcal F$ turns out to be a 
solution to the problem \pier{made precise by \eqref{pier2}--\eqref{eq:22}}. To this aim, we are going to show that some power $\mathcal F^j$ ($j\in{\mathbb N}$) is a contract\pier{ion mapping} in $S$ and, as a consequence, it admits a unique fixed point, which results at the end to provide the unique solution to our system.\finele

To this \pier{aim}, we pick a pair \pier{$\{v_i\}_{i=1,2}\subset S$ and set}
\begin{equation*}
  u_i:=\mathcal F(v_i), \quad \xi_i := \mu_i -\mu_\flat-\psi_*'(v_i) -\partial_t u_i, 
\end{equation*}
where $\mu_i $ is the solution to \eqref{punto2} corresponding to $v_i$, $i=1,2$. 
\pier{Then, it is easy to verify that
\begin{subequations}\label{pier11}
\begin{gather}
\pier{\partial_t u_i }(t)+A\mu_i(t)=0\quad \textrm{in }V',\, \hbox{ for all }t\in[0,T],\label{pier11a}\\
\pcol{ \mu_i (t) = (\pier{\partial_t u_i }+\xi_i+\mu_\flat+\psi'(v_i))(t)\quad \hbox{a.e. in } \Omega , \, \hbox{ for all }t\in[0,T],}, \label{pier11b}\\
\pcol{\xi_i (t) \in\beta(\pier{\partial_t u_i }(t) )\quad \hbox{a.e. in } \Omega , \, \hbox{ for all }t\in[0,T],}\label{pier11c}\\
u_i (0)=u_0\quad\textrm{a.e. in } \Omega \label{pier11d}
\end{gather}
\end{subequations}
for $i=1,2$. Now, we use the notation $\bar u$ for the difference of $u_1-u_2$, and the same notation for $\bar\xi$, $\bar\mu$ \pier{and $\bar v$.} We take the difference of \giu\eqref{pier11a} for $i=1,2$, \finegiu 
test it by $A^{-1} (\partial_t\bar u \pcol{(t)} ) $ and, at the same time, we test the difference of \eqref{pier11b} by $\partial_t \bar u \pcol{(t)}$. Then we combine the obtained equalities and use the properties of \giu $A^{-1}$ stated \finegiu in \eqref{eq:20bis} and  \eqref{eq:20tris}.
Hence, we have that}\elena
\begin{align}\label{eq:39}
& \|\pier{\partial_t {\bar u}}(t)\|_{*}^2+\|\pier{\partial_t {\bar u}}(t)\|_{H}^2+\int_\Omega\bar\xi(t)\pier{\partial_t {\bar u}}(t)\le \int_\Omega |\psi_*'(v_1(t))-\psi_*'(v_2(t))||\pier{\partial_t {\bar u}(t)}|\nonumber\\
&\le \frac 12\|\bar u_t(t)\|^2_H+\frac 12\int_\Omega|\psi_*'(v_1(t))-\psi_*'(v_2(t))|^2.
\end{align}
Due to \eqref{eq:52}, we handle the right hand \pcol{side} observing that
\begin{equation}\label{eq:45}
\int_\Omega |\psi_*'(v_1(t))-\psi_*'(v_2(t))|^2\le L^2\!\int_\Omega\pier{|\bar v(t)|}^2.
\end{equation}
In addition, by \pier{the monotonicity of $\beta$ and \eqref{pier11c} we deduce that}
\begin{equation*}
\int_\Omega\bar\xi(t)\pier{\partial_t {\bar u}}(t)\geq0.
\end{equation*}
Thus, we \pier{easily obtain}
\begin{equation}
\|\pier{\partial_t {\bar u}}(t)\|_H\le \pier{L} \|\bar v(t)\|_H
\end{equation}
and, \pier{as both $u_1$ and $u_2$ satisfy the same initial condition~\eqref{pier11d},  
we can easily infer that
\begin{equation}\label{eq:51prima}
\|\mathcal F(v_1(t))-\mathcal F(v_2(t))\|_H\le \pier{L} \int_0^t\|v_1(s)-v_2(s)\|_H\pier{{\rm d}s} \quad \pier{\hbox{for all } \, t\in [0, T].} 
\end{equation}
This inequality leads to}
\begin{equation}\label{eq:51}
   \| \mathcal F(v_1)-\mathcal F(v_2)\|_{C^0([0,t];H)}\le \pier{L}  t\|v_1-v_2\|_{C^0([0,t];H)} \quad \pier{\hbox{for all } \, t\in [0, T].}
   \end{equation}
An \pier{iteraction} of the argument, due to \eqref{eq:51prima} and \eqref{eq:51}, \finele \pier{leads} to
\elena
\begin{align*}
  \|\mathcal F^2(v_1)-\mathcal F^2(v_2)\|_{C^0([0,t];H)}\le \pier{L} 
  \int_0^t \pier{\|\mathcal F(v_1)(s)-\mathcal F(v_2)(s)\|_{H}}\pier{{\rm d}s}
\\
\le \pier{L} ^2\int_0^t s\|v_1-v_2\|_{C^0([0,s];H)}\pier{{\rm d}s}\nonumber
\pier{{}\le\frac {(\pier{L} t)^2}{2}\|v_1-v_2\|_{C^0([0,t];H)}}.
\end{align*}
\finele
By iterating $j$ times, we find $\|\mathcal F^j(v_1)-\mathcal F^j(v_2)\|_{C^0([0,t];H)}\le\frac {(\pier{L} t)^j}{j !}\|v_1-v_2\|_{C^0([0,t];H)}$ for all $t\in [0,T]$, whence, in particular,
\begin{align}
  \|\mathcal F^j(v_1)-\mathcal F^j(v_2)\|_{C^0([0,T];H)}&\le\frac {(\pier{L} T)^j}{j!}\|v_1-v_2\|_{C^0([0,T];H)}.
\end{align}
Thus, for $j$ large enough $\mathcal F^j$ \elena turns out to be \finele a contraction mapping from $S$ into itself, as announced; as a consequence, $\mathcal F^j$ has a \pier{unique} fixed point $u^*$, which is also \pier{the unique} fixed point for $\mathcal F$. \pier{In view of \eqref{pier11}, this} fixed point yields the triplet $(u^*,\xi^*,\mu^*)$ that solves \elena the problem \pier{\eqref{eq:22} in which} $\psi$ is substituted by $\psi_*$. \finele
Of course, \pier{for $(u^*,\xi^*,\mu^*)$} we can repeat the estimates carried out in Section~\ref{sec:basic-estimates}. In particular --- and this is the crucial point --- we can derive for $\mu^*$ the same estimate as \eqref{eq:49}, namely,
\begin{equation}\label{eq:68}
  \|\mu^*-\mu_\flat\|_{L^\infty(Q)}\le M,
\end{equation}
with the \emph{same value of the constant} $M$. In fact, \pier{if one checks carefully the estimates, \giu one can \finegiu see that $\psi_*$ appears in \eqref{eq:46} with the integral of 
$\psi_* (u_0) \equiv \psi (u_0)$ and in \eqref{eq:41} with $\psi'_* (u_0) \equiv \psi' (u_0)$ (cf.~\eqref{eq:6} and \eqref{eq:33a}), and especially with the constant $K_1$ in \eqref{eq:23c} which, as observed at the beginning of this section, can be the same for $\psi$ and $\psi_*$.} 
\elena In addition, by \finele
\gius{its very definition \pier{the derivative $\psi'_*$} satisfies}
\begin{subequations}\label{eq:48}
\begin{gather}
  \psi'_*(r)\ge M\quad\textrm{for all }r\ge k^*,\label{eq:48a}\\ 
  \psi'_*(r)\le -M\quad\textrm{for all }r\le k_*,\label{eq:48b}
\end{gather}
\end{subequations} 
\pier{as $\psi'$ does in \eqref{eq:33b}--\eqref{eq:33c}. Thus, a repetition of the argument leading to \eqref{eq:28} yields
\begin{equation}\label{eq:53}
  k_*\le u^*\le k^*\quad\textrm{a.e. in }Q
\end{equation}
and, since $\psi=\psi_*$ in $[k_*,k^*]$, \eqref{eq:53} entails
\begin{equation}
  \psi_*'(u^*)=\psi'(u^*)\quad\textrm{ a.e. in }Q.
\end{equation}
In other words,  \elena$(u^*,\xi^*,\mu^*)$ \finele  is \pier{actually} a solution to the original problem \eqref{eq:22} and fulfills~\eqref{pier10}. Moreover, $(u^*,\xi^*,\mu^*)$ is the unique solution of \eqref{eq:22},
owing to the continuous dependence property~\eqref{contdep} 
proved in Section~\ref{sec:cont-depend-data-1}. Finally, recalling the smoothness of  
$\Omega$ and $\Gamma$ and the homogeneous boundary condition  on $\Gamma$, by \eqref{pier2}, \eqref{eq:6-2} and standard 
elliptic regularity estimates we obtain \eqref{pier1} and the regularity  $\mu \in C^0([0,T];W)$ for  $(u^*,\xi^*,\mu^*)$. Therefore, Theorem~\ref{Th1} is completely proved.}

\pier{\begin{remark}
\label{rem-pier}
{\rm Note that if we assume $u_0 \in H^1(\Omega) $ and $ \mu_\flat\in L^2(0,T;H^1(\Omega))$
besides \eqref{eq:6} and \eqref{eq:24}, then we can recover the additional regularity 
$u\in H^1(0,T; H^1(\Omega))$ for the solution component $u$. Indeed, it suffices to take formally the gradient of \eqref{pier8} and test it by $\nabla (\partial_t u)$. Thanks to the Lipschitz continuity (with Lipschitz constant $1$) of $(I+\beta)^{-1}$ and of $\psi'_*$ (which can replace $  \psi'$ as we have seen) with constant $L$ (cf.~\eqref{eq:52}), by the Young inequality we easily infer that
\begin{align*} 
\frac12 \Vert \nabla (\partial_t u) (t) \Vert_H^2 \leq 
\Vert \nabla (\mu-\mu_\flat) (t) \Vert_H^2 + L^2 \Vert \nabla  u (t) \Vert_H^2.
\end{align*} 
Now, pointing out that  $\mu-\mu_\flat \in L^2(0,T;H^1(\Omega))$ (cf., e.g., \eqref{eq:6-1}), as 
$$\| \nabla u(t)\|_H^2\le 2 \| \nabla u_0\|^2_{H}+ 2T \!\int_0^t\|\nabla (\pier{\partial_t {u}})(s)\|_H^2\,\pier{{\rm d}s}, $$ we can easily apply the Gronwall lemma and find out that  $\partial_t u \in L^2(0,T; H^1(\Omega))$.} 
\end{remark}}
\def\cprime{$'$}


\end{document}